\documentclass{gtmon_a}
\pdfoutput=1

\usepackage[T1,T5]{fontenc}

\proceedingstitle{Proceedings of the School and Conference in Algebraic
Topology (The Vietnam National University, Hanoi, 9-20 August 2004)}
\conferencestart{9 August 2004}
\conferenceend{20 August 2004}
\conferencename{School and Conference in Algebraic Topology}
\conferencelocation{Vietnam National University, Hanoi, Vietnam}

\editor{John Hubbuck}
\givenname{John}
\surname{Hubbuck}

\editor{Nguy\~\ecircumflex{}n H\,V H\uhorn{}ng}
\givenname{Nguy\~\ecircumflex{}n H V}
\surname{H\uhorn{}ng}

\editor{Lionel Schwartz}
\givenname{Lionel}
\surname{Schwartz}

\title[The action of the primitive Steenrod--Milnor operations]{The action of the primitive Steenrod--Milnor\\ operations on the modular invariants}

\author{{\fontencoding{T5}\selectfont Nguy\~\ecircumflex{}n Sum}}
\surname{Sum}
\givenname{Nguy\~\ecircumflex{}n}
\address{Department of Mathematics\\
University of Quynhon\\\newline
170 An Duong Vuong\\Quynhon, Binhdinh\\Vietnam}
\email{ngnsum@yahoo.com}
\urladdr{}

\volumenumber{11}
\issuenumber{}
\publicationyear{2007}
\papernumber{16}
\startpage{349}
\endpage{367}

\doi{}
\MR{}
\Zbl{}

\arxivreference{}


\received{15 November 2004}
\revised{28 August 2005}
\accepted{12 December 2005}
\published{14 November 2007}
\publishedonline{14 November 2007}
\proposed{}
\seconded{}
\corresponding{}
\editor{}
\version{}

\makeatletter
\def\cnewtheorem#1[#2]#3{\newtheorem{#1}{#3}[section]
\expandafter\let\csname c@#1\endcsname\c@thm}

\AtBeginDocument{\let\bar\wbar\let\tilde\wtilde\let\hat\what}
\makeop{St}

\theoremstyle{plain}
\newtheorem{thm}{Theorem}[section]
\cnewtheorem{lem}[thm]{Lemma}
\cnewtheorem{corl}[thm]{Corollary}
\cnewtheorem{prop}[thm]{Proposition}
\theoremstyle{definition}
\cnewtheorem{defn}[thm]{Definition}
\cnewtheorem{rem}[thm]{Remark}

\makeatother  

\begin{document} 

\begin{htmlabstract}
We compute the action of the primitive Steenrod&ndash;Milnor operations
on generators of algebras of  invariants of subgroups of general linear
group GL<sub>n</sub>=GL(n,F<sub>p</sub>) in the polynomial algebra with
p an odd prime number.
\end{htmlabstract}

\begin{abstract}
We  compute the action of the primitive Steenrod--Milnor
operations on generators of algebras of  invariants of subgroups of
general linear group ${GL_n=GL(n, \mathbb F_p)}$ in the polynomial
algebra with $ p$ an odd prime number.
\end{abstract}

\maketitle

\section{Introduction}
\label{sec1}
Let  $p$ be an odd prime number. Denote by $GL_n  =  GL(n,\mathbb F_p)$ the general linear  group over the prime field $\mathbb F_p$ and $T_n$ the Sylow $p$--subgroup of $GL_n$ consisting of all upper triangular matrices with 1 on the main diagonal. For any integer  $d,\ 1 \le
d \le p-1$, we set
$$SL_n^d = \{ \omega \ \in GL_n;\ (\det \omega  )^d = 1 \}.$$
It is easy to see that $SL_n^d$ is a subgroup of $GL_n$ and $SL_n^{p-1} = GL_n$. Each subgroup of $GL_n$ acts on  $P_n=E(x_1,\ldots,x_n)\otimes \mathbb F_p(y_1,\ldots,y_n)$ in the usual manner. Here and in what follows, $E(.,\ldots,. )$ and $\mathbb F_p(.,\ldots,.)$ are the exterior and polynomial algebras over $\mathbb F_p$ generated by the indicated variables. We grade $P_n$ by assigning $\dim x_i=1$ and $\dim y_i=2.$

Dickson \cite{1} showed that the invariant algebra $\smash{ \mathbb
F_p(y_1,\ldots,y_n)^{\text{\tiny{$GL_n$}}}}$ is a polynomial algebra
generated by the Dickson invariants $Q_{n,s},\ 0\le s<n$. Hu\`ynh M\`ui
\cite{2,3} computed the invariant algebras
$\smash{P_n^{\text{\tiny{$T_n$}}}}$  and
$\smash{P_n^{\text{\tiny{$SL_n^d$}}}}$ for $d=1, p-1, (p-1)/2$. He proved
that $\smash{P_n^{\text{\tiny{$T_n$}}}}$ is generated by $V_m$, $1 \le m
\le n$,  $M_{m, s_1,  \ldots, s_k}$, $0 \le s_1 < \ldots < s_k < m\le n$
and that $\smash{P_n^{\text{\tiny{$SL_n^d$}}}}$ is generated by  $L_n^d,\
Q_{n, s},\ 1 \le  s <  n,\ \smash{M_{n, s_1,  \ldots, s_k}^{(d)}},\  0 \le
s_1 < \ldots < s_k < n.$ Here $V_m,  \smash{M_{n, s_1,  \ldots,
s_k}^{(d)}}$ are M\`ui invariants and $L_n^d, Q_{n,s}$ are  Dickson invariants (see \hbox{\fullref{sec2}}). Note that $\smash{M_{n, s_1,  \ldots, s_k}^{(1)}}=M_{n, s_1,  \ldots, s_k}$.

Let $\mathcal{A}(p)$ be the mod $p$ Steenrod algebra and let
$\tau_s, \xi_i$ be the Milnor elements of dimensions $2p^s-1,\
2p^i-2$ respectively in the dual algebra $\mathcal{A}(p)^*$ of $\mathcal{A}(p)$.
In \cite{7}, Milnor showed that as an algebra,
$$\mathcal{A}(p)^* = E(\tau_0,\tau_1,\ldots )\ \otimes\ \mathbb F_p(\xi_1,\xi_2,\ldots ). $$
Then $\mathcal{A}(p)^*$ has a basis consisting of all monomials
    $$\tau_S\xi^R \ =\ \tau_{s_1}\ldots \tau_{s_k}\xi_1^{r_1}\ldots
     \xi_m^{r_m},$$
with $S = (s_1,\ldots ,s_k),\ 0 \le s_1 <\ldots <s_k ,
 R = (r_1,\ldots ,r_m),\ r_i \ge 0 $. Let $\St^{S,R} \in \mathcal{A}(p)$
denote the dual of $\tau_S\xi^R$ with respect to that basis.
Then $\mathcal{A}(p)$ has a basis consisting of all operations $\St^{S,R}$. For $S=\emptyset, R=(r)$, $\St^{\emptyset, (r)}$ is nothing but the Steenrod operation $P^r$. So, we call $\St^{S,R}$ the Steenrod--Milnor operation of type $(S,R)$.

We have the Cartan formula
$$\St^{S,R}(uv) =\sum_{\overset{\scriptstyle{S_1\cup S_2=S}}{R_1+R_2=R}}(-1)^{(\dim u+\ell(S_1))\ell(S_2)}(S:S_1,S_2)\St^{S_1,R_1}(u)\St^{S_2,R_2}(v),$$
where $ R_1=(r_{1i}),\ R_2=(r_{2i}),\ R_1+R_2=(r_{1i}+r_{2i}), S_1\cap S_2=\emptyset, u,v\in P_n$, $\ell(S_i)$ is the length of $S_i$ and
$$(S:S_1,S_2)= \text{sign}\begin{pmatrix} s_1&\ldots &s_h& s_{h+1}&\ldots &s_k\\
s_{11}&\ldots &s_{1h}&s_{21}&\ldots &s_{2r}\end{pmatrix},$$ 
with $S_1=(s_{11},\ldots,s_{1h}), s_{11}<\ldots <s_{1h}$, $S_2=(s_{21},\ldots,s_{2r}), s_{21} < \ldots < s_{2r}$ (see Mui \cite{3}).

We denote 
$ \St_u=\St^{(u),(0)},\ \St^{\Delta_i}=\St^{\emptyset,\Delta_i}$, where $\Delta_i=(0,\ldots,1,\ldots,0)$ with 1 at the $i$--th place.  In  \cite{3}, Hu\`ynh M\`ui proved that as a coalgebra,
$$\mathcal{A}(p) =\Lambda(\St_0,\St_1,\ldots)\otimes \Gamma(\St^{\Delta_1},\St^{\Delta_1},\ldots ).$$
Here, $\Lambda(\St_0,\St_1,\ldots)$ (resp. $\Gamma(\St^{\Delta_1},\St^{\Delta_2},\ldots )$) denotes the exterior (resp.\ polynomial) Hopf algebra with divided powers generated by the  primitive Steenrod--Milnor operations $\St_0,\ \!\St_1,\ldots$ (resp. $\St^{\Delta_1}, \St^{\Delta_2},\ldots )$. 

The Steenrod algebra  $\mathcal{A}(p)$ acts on $P_n$ by means of the Cartan formula together with the relations 
\begin{align*} 
&\St^{S,R}x_k=\begin{cases} x_k, & S=\emptyset,\ R=(0),\\
y_k^{p^u}, &S=(u),\ R=(0), \\ 
0,  &\text{otherwise}, 
\end{cases}\tag{\hbox{i}}\label{(i)}\\
&\St^{S,R}y_k=\begin{cases}
 y_k, & S=\emptyset ,\ R=(0),\\
y_k^{p^i},&S=\emptyset,\ R=\Delta_i, \\ 
0, &\text{otherwise},  \end{cases}\tag{\hbox{ii}}\label{(ii)}
\end{align*}
for $k=1,\ \! 2, \ldots, n$ (see Steenrod and Epstein \cite{13} and Sum \cite{11}). 
Since this action commutes with  the action of $GL_n$, it induces  actions of $\mathcal{A}(p)$ on $\smash{P_n^{\text{\tiny{$T_n$}}}}$ and $\smash{P_n^{\text{\tiny{$SL_n^d$}}}}$. 

The action of $\St^{S,R}$ on the modular invariants of subgroups of general linear group has been studied by many authors. This action for $S=\emptyset,\ \! R=(r)$ was explicitly determined by H\uhorn{}ng and Minh \cite{5}, Kechagias \cite{4}, Madsen and Milgram \cite{6} and \hbox{Sum \cite{11}}. Smith and Switzer \cite{12}, Wilkerson \cite{14} and Neusel \cite{8} have studied the action of $\St^{\Delta_i}$ on the Dickson invariants.

The purpose of the paper is to compute the action of the primitive Steenrod--Milnor operations on generators of $\smash{P_n^{\text{\tiny{$T_n$}}}}$ and $\smash{P_n^{\text{\tiny{$SL_n^d$}}}}$.

The rest of the paper contains three sections. In \fullref{sec2}, we recall some needed information on the invariant theory and compute the action of the primitive Steenrod--Milnor operations on the determinant invariants. In \fullref{sec3}, we compute the action of  the primitive Steenrod--Milnor operations on Dickson and M\`ui invariants. Finally, we give in \fullref{sec4} some formulae from which we can describe our results in terms of Dickson and M\`ui invariants.

{\bf Acknowledgements}\qua  
 The author is grateful to the referee for his valuable comments on the first manuscript of this paper.

\section{Preliminaries}
\label{sec2}

\begin{defn} 
Let $(e_{k+1 },\ldots,e_m),\ 0  \leq   k < m  \leq  n$, be a sequence of nonnegative integers.  Following Dickson \cite{1} and M\`ui \cite{2}, we define  
$$ [k;e_{k+1},  \ldots,  e_m]  =  \frac 1{k!}
\vmatrix x_1&\cdots &x_m\\
  \vdots&\cdots &\vdots\\
  x_1&\cdots &x_m\\
  y_1^{p^{e_{k+1}}}&\cdots &y_m^{p^{e_{k+1}}}\\
  \vdots&\cdots  &\vdots\\
   y_1^{p^{e_m}} & \cdots & y_m^{p^{e_m}}
   \endvmatrix. $$
\end{defn}
 The precise meaning of the right hand side is given in \cite{2}. For $k=0$, we write
$$[0;e_{1},  \ldots,  e_m] =  [e_{1},  \ldots,  e_m] =\det(y_i^{p^{e_j}}).$$
In particular, we set
 \begin{align*} L_{m,s}&=[0,1,  \ldots,\hat s,\ldots,  m],\ \! 0\le s \le m \le n,\\ 
L_m =L_{m,m}&=[0,1,\ldots,m-1],\ \  
\\ 
M_{m,s_1,\ldots,s_k}  &= [k;0,\ldots,\hat s_1,\ldots,\hat
s_k, \ldots, m -1], \end{align*}
for $0 \leq \ s_1 < \ldots < s_k < m \le n$. Each $[k;e_{k+1},  \ldots,  e_m]$ is an invariant of $SL_m^1$ and $[e_{1},  \ldots,  e_m]$ is divisible by $L_m$. Then, Dickson invariants $Q_{n,s}$ and M\`ui invariants $ M_{n,s_1,\ldots,s_k}^{(d)}$ and $V_m$ are defined by
$$ Q_{n,s}=L_{n,s}/L_n,\quad M_{n,s_1,\ldots,s_k}^{(d)} = M_{n,s_1,\ldots,s_k}L_n^{d-1} \quad\text{and}\quad V_m = L_m/L_{m-1}.$$
Here, by convention, $L_0=[\emptyset] =1.$

Now we prepare some data in order to prove our main results.

\begin{lem}\label{lem2.2} Suppose $e_\ell\ne e_j$ for $\ell\ne j$, $u \ge 0$. Then we have
$$\St_u[k;e_{k+1},\ldots,e_n] = \begin{cases} (-1)^{k-1}[k-1;u,e_{k+1},\ldots,e_n],&k>0,\\
0, & k=0.\end{cases}$$
\end{lem}

\begin{proof} Let $I$ be a subset of $\{1,\ldots,n\}$ and let $I'$ be its complement
in $\{1,2,\ldots,n\}$. 
Writing 
$I =  \{i_1,i_2,\ldots,i_k\}$ and $I' =  \{i_{k+1},i_{k+2},\ldots,i_n\}$ with $i_1 < i_2 < \ldots < i_k$ and $i_{k+1} < i_{k+2} < \ldots < i_n$. We set 
\begin{align*} x_I &= x_{i_1}x_{i_2}\ldots x_{i_k}, \\ {[e_{k+1},e_{k+2},\ldots, e_n ]_I} &=  {[ e_{k+1},e_{k+2},\ldots, e_n](y_{i_{k+1}},y_{i_{k+2}},\ldots, y_{i_n})}\\
\sigma_I&= \begin{pmatrix} 1& 2& \ldots & n\\ i_1& i_2 & \ldots & i_n\end{pmatrix} \in \Sigma_n, \end{align*} where $\Sigma_n$ is the symmetric group on $n$ letters. Using the Laplace development, we have
$$ [k;e_{k+1},e_{k+2},\ldots,e_n] = \sum_I\text{sign}\  \sigma_I x_I[e_{k+1},e_{k+2},\ldots,e_n]_I.$$
From the relation \ref{(ii)}, we see that $\St_u[e_{k+1},e_{k+2},\ldots,e_n]_I =0$. Then, using the Cartan formula, we get
\begin{equation}\label{ct3} \St_u[k;e_{k+1},e_{k+2},\ldots,e_n] = \sum_I\text{sign}\ \sigma_I \St_u(x_I)[e_{k+1},e_{k+2},\ldots,e_n]_I.\end{equation}
In \cite[5.2]{3}, M\`ui showed that $$\St_u(x_I) = (-1)^{k-1}[k-1;u](x_{i_1},x_{i_2},\ldots ,x_{i_k},y_{i_1},y_{i_2},\ldots ,y_{i_k}).$$ Hence, using \eqref{ct3} and the Laplace development we obtain the lemma.  
\end{proof}

\begin{lem}\label{lem2.3} Suppose $e_\ell\ne e_j$ for $\ell \ne j$, $e_{k+1}<e_j$ for $j>k+1$. Then we have
$$\St^{\Delta_i}[k;e_{k+1},\ldots,e_n] =\begin{cases} [k;i,e_{k+2},\ldots,e_n], &e_{k+1}=0,\\ 0, &e_{k+1}>0.\end{cases}$$
\end{lem}

\begin{proof} Suppose $e_{k+1}>0$. From the relations \ref{(i)} and \ref{(ii)} and the Cartan formula, we easily obtain 
$$ \St^{\Delta_i}x_\ell=0,\ \ \St^{\Delta_i}y_\ell^{p^{e_j}}= p^{e_j}y_\ell^{p^{e_j}+p^i-1}=0,$$ 
for  $\ell = 1, 2,\ldots , n$ and $j = k+1,k+2,\ldots, n.$
 From this, we get $$\St^{\Delta_i}[k;e_{k+1},\ldots, e_n] = 0.$$
If $e_{k+1}=0$ then $\St^{\Delta_i}y_\ell^{p^{e_j}}=0,$ for $ \ell = 1, 2,\ldots , n$ and $j = k+2,\ldots, n$, and $$\St^{\Delta_i}y_\ell^{p^{e_{k+1}}} = \St^{\Delta_i}y_\ell = y_\ell^{p^i}.$$ Hence, using the Laplace development and the Cartan formula, we obtain
$$\St^{\Delta_i}[k;e_{k+1},e_{k+2},\ldots, e_n] = [k;i,e_{k+2},\ldots, e_n].\proved$$
\end{proof}

To make the paper self-contained, we give here a proof for the following theorem, which will be needed in the next section.

\begin{thm} \label{thm2.4}{\rm (Sum \cite{10})}\qua  Let $(e_{1},e_{2},\ldots, e_n)$ be a sequence of nonnegative integers and $0\le k < n$. We have
\begin{align}&\notag [e_{1},e_{2},\ldots,e_{n-1}, e_n+n-1] \\ 
&\label{ct5}\ \  = \sum_{s=0}^{n-2}(-1)^{n+s}[e_{1},e_{2},\ldots,e_{n-1}, e_n+s]Q_{n-1,s}^{p^{e_n}}
+ [e_{1},e_{2},\ldots, e_{n-1}]V_n^{p^{e_n}},\\
&\label{ct6}[k;e_{k+1},\ldots , e_{n-1}, e_n+n] = \sum_{s=0}^{n-1}(-1)^{n+s-1}[k;e_{k+1},\ldots,e_{n-1}, e_n+s]Q_{n,s}^{p^{e_n}}.
\end{align}
\end{thm}

\begin{proof} We recall M\`ui's formula in \cite{2},
\begin{multline*}[k;e_{k+1},\ldots, e_n] =\\ (-1)^{k(k-1)/2}\sum_{0\le s_1<\ldots <s_k}(-1)^{s_1+\ldots +s_k}M_{n,s_1, \ldots,s_k}[s_1, \ldots, s_k, e_{k+1} ,\ldots, e_n]/L_n.
\end{multline*}
Hence, it suffices to prove the theorem for $k=0$.

The proof of the theorem proceeds by induction on $n$. It is easy to see that the theorem holds for $n=1$. Suppose $n\ge 2$ and the theorem holds for $n-1$.

Using the Laplace development and the inductive hypothesis, we have
\begin{align*}
[e_{1},&\ldots,e_{n-1}, e_n+n-1]\\ &= \sum_{t=1}^{n-1}(-1)^{n+t}[e_{1},\ldots,\hat e_t,\ldots, e_{n-1}, e_n+n-1] y_n^{p^{e_t}} + [e_1,\ldots,e_{n-1}]y_n^{p^{e_n+n-1}}  \\
& = \sum_{t=1}^{n-1}(-1)^{n+t}\Big(\sum_{s=0}^{n-2}(-1)^{n+s}[e_{1},\ldots,\hat e_t,\ldots, e_{n-1}, e_n+s]Q_{n-1,s}^{p^{e_n}}\Big) y_n^{p^{e_t}}\\
&\hskip234pt + [e_1,\ldots,e_{n-1}]y_n^{p^{e_n+n-1}}\\
&=\sum_{s=0}^{n-2}(-1)^{n+s}\Big( \sum_{t=1}^{n-1}(-1)^{n+t}[e_{1},\ldots,\hat e_t,\ldots, e_{n-1}, e_n+s]y_n^{p^{e_t}}\Big)Q_{n-1,s}^{p^{e_n}} \\
&\hskip234pt + [e_1,\ldots,e_{n-1}]y_n^{p^{e_n+n-1}}\\
& =\sum_{s=0}^{n-2}(-1)^{n+s}[e_{1},\ldots,e_{n-1}, e_n+s]Q_{n-1,s}^{p^{e_n}}\\
&\hskip4.95cm  + [e_1,\ldots,e_{n-1}]\sum_{s=0}^{n-1}(-1)^{n+s-1}Q_{n-1,s}^{p^{e_n}}y_n^{p^{e_n+s}}.
\end{align*}
Since $V_n=\sum_{s=0}^{n-1}(-1)^{n+s-1}Q_{n-1,s}y_n^{p^{s}}$, the relation \eqref{ct5} holds for $n$.

Now we prove the relation \eqref{ct6} for $n$. By a direct calculation using \eqref{ct5} and the relation $Q_{n,s}=Q_{n-1,s-1}^p + Q_{n-1,s}V_n^{p-1}$, we get
 \begin{align*}[e&_{1},e_2,\ldots , e_{n-1}, e_n+n]\\
 &= \sum_{s=1}^{n-1}(-1)^{n+s-1}[e_{1},\ldots,e_{n-1}, e_n+s]Q_{n-1,s-1}^{p^{e_n+1}}+ [e_1,\ldots,e_{n-1}]V_n^{p^{e_n+1}}\\
 &= \sum_{s=1}^{n-1}(-1)^{n+s-1}[e_{1},\ldots,e_{n-1}, e_n+s]Q_{n,s}^{p^{e_n}}\\
&\hskip5cm - [e_1,\ldots,e_{n-1},e_n+n-1]V_n^{(p-1)p^{e_n}}\\
&\quad + \Big(\sum_{s=1}^{n-2}(-1)^{n+s}[e_1,\ldots,e_{n-1},e_n+s]Q_{n-1,s}^{p^{e_n}} 
+ [e_1,\ldots,e_{n-1}]V_n^{p^{e_n}}\Big)V_n^{(p-1)p^{e_n}}.
\end{align*}
Combining this equality and the relation \eqref{ct5} we obtain
\begin{align*}
[e_{1},e_2,\ldots , e_{n-1}, e_n+n]
 = \sum_{s=1}^{n-1}&(-1)^{n+s-1}[e_{1},\ldots,e_{n-1}, e_n+s]Q_{n,s}^{p^{e_n}}\\
- &(-1)^{n}[e_1,\ldots,e_{n-1},e_n]Q_{n-1,0}^{p^{e_n}} V_n^{(p-1)p^{e_n}}.
\end{align*}
Since $Q_{n,0}=Q_{n-1,0}V_n^{p-1}$, the relation \eqref{ct6} holds for $n$.

This completes the proof of \fullref{thm2.4}.   
\end{proof}

\section{Main results}
\label{sec3}
Observe that using the Cartan formula and the relations \ref{(i)} and \ref{(ii)}, we obtain $\St_ux=0$ for either $x=Q_{n,s}$ or $x=V_n$. So, in this section we only compute $\St^{\Delta_i}x$ for $x= Q_{n,s}, V_n, \smash{M_{n,s_1,\ldots,s_k}^{(d)}}$ and $\smash{\St_uM_{n,s_1,\ldots,s_k}^{(d)}}$.

\begin{thm}\label{thm3.1} For any integers $i,\ \! n, \ \! s$ with $0\le s <n$ and $i \ge 1$, we have
$$\St^{\Delta_i}Q_{n,s}= (-1)^n[0,1,\ldots,\hat s, \ldots,n-1,i]L_n^{p-2}.$$
\end{thm}

\begin{proof} Since $L_{n,s}=L_nQ_{n,s}$, using the Cartan formula, we get 
\begin{equation} \label{ct4} \St^{\Delta_i}L_{n,s}= L_n\St^{\Delta_i}Q_{n,s}+Q_{n,s}\St^{\Delta_i}L_{n}.
\end{equation}
According to \fullref{lem2.3}, we have
$$\St^{\Delta_i}L_{n,s}= \begin{cases} [i,1,2,\ldots,\hat s,\ldots,n],&s>0,\\ 0,&s=0.\end{cases}$$
In particular,  $\St^{\Delta_i}L_{n}=[i,1,2,\ldots,n-1]$. 

 If  $s=0$ then $\St^{\Delta_i}L_{n,s}=0$ and 
\begin{align*}\St^{\Delta_i}L_{n}&=[i,1,2,\ldots,n-1]\\ 
& = (-1)^{n-1}[1,2,\ldots,n-1,i].
\end{align*}
Combining \eqref{ct4} and the above equalities,  we get
\begin{align*} \St^{\Delta_i}Q_{n,0}&= -(\St^{\Delta_i}L_n) Q_{n,0}/L_n\\
&=(-1)^{n}[1,2,\ldots,n-1,i]Q_{n,0}/L_n.
\end{align*}
Since $Q_{n,0}=L_n^{p-1}$, the theorem holds.

 If  $s>0$ then $\St^{\Delta_i}L_n=[i,1,2,\ldots,n-1]$ and 
$\St^{\Delta_i}L_{n,s} = [i,1,2,\ldots,\hat s,\ldots, n].$
 Hence, using \fullref{thm2.4}, we get
\begin{align*}\St^{\Delta_i}L_{n,s}  &=\sum_{t=0}^{n-1}(-1)^{n-1+t}[i,1,2,\ldots,\hat s,\ldots,n-1, t]Q_{n,t}\\
&=(-1)^{n-1} [i,1,2,\ldots,\hat s,\ldots,n-1, 0]Q_{n,0}\\
&\qquad +(-1)^{n-1+s}[i,1,2,\ldots ,\hat s,\ldots ,n-1, s]Q_{n,s}\\
&= [i,1,2,\ldots,n-1]Q_{n,s}-[i,0,1,\ldots,\hat s,\ldots,n-1]Q_{n,0}.
\end{align*}
Combining \eqref{ct4}, the above equalities and  the relation $Q_{n,0}=L_n^{p-1}$, we get
$$\eqalignbot{\St^{\Delta_i}Q_{n,s}&= (\St^{\Delta_i}L_{n,s}-Q_{n,s}\St^{\Delta_i}L_{n})/L_n\cr
&=-[i,0,1,2,\ldots,\hat s,\ldots,n-1]Q_{n,0}/L_n\cr
&=(-1)^n[0,1,2,\ldots,\hat s,\ldots, n-1,i]L_n^{p-2}.}
\proved$$
\end{proof}

The following  was proved in Smith and Switzer \cite{12}  by another method.

\begin{corl}{\rm (Smith--Switzer \cite{12})}\qua For any integers $i,\ \! n, \ \! s$ with $0\le s <n$ and $1\le i\le n$, we have
$$\St^{\Delta_i}Q_{n,s}=\begin{cases} (-1)^{s-1}Q_{n,0},& i=s>0,\\
(-1)^{n}Q_{n,0}Q_{n,s},&i=n,\\
0,\qquad &\text{otherwise.}\end{cases}$$
\end{corl}

\begin{proof} Suppose  $i = s$. According to \fullref{thm3.1}, we have 
\begin{align*}
\St^{\Delta_s}Q_{n,s}&=(-1)^n[0,1,\ldots ,\hat s,\ldots, n-1,s]L_n^{p-2}\\
&=(-1)^{s-1}[0,1,\ldots,n-1]L_n^{p-2}\\ 
&= (-1)^{s-1}L_n^{p-1} = (-1)^{s-1}Q_{n,0}.
\end{align*}
If $i<n$ and $i\ne s$ then $[0,1,\ldots,\hat s,\ldots,n-1,i]=0$. Hence, $\St^{\Delta_i}Q_{n,s} = 0.$
\begin{align*}
\tag*{\hbox{If $i=n$ then}}\St^{\Delta_n}Q_{n,s} &= (-1)^n[0,1,\ldots,\hat s,\ldots,n-1,n]L_n^{p-2}\qquad\\
&=(-1)^nL_{n,s}L_n^{p-2}\\ 
& = (-1)^nL_n^{p-1}Q_{n,s}\\ 
& = (-1)^nQ_{n,0}Q_{n,s}.
\end{align*}
The corollary follows.
\end{proof}

Now, we show that  our formula in \fullref{thm3.1} implies Wilkerson's formula  in \cite{14}. To do this, we need the following.

\begin{prop}{\rm (Sum \cite{10})}\qua\label{prop3.3}  Let $(e_{k+1},e_{k+2},\ldots, e_n)$ be a sequence of nonnegative integers with $0\le k <n$ and $e_\ell \ne e_j$ for $\ell \ne j$. Then
$$P^r[k;e_{k+1},e_{k+2},\ldots , e_n] = \begin{cases} [k;e_{k+1}+\varepsilon_{k+1},e_{k+2}+\varepsilon_{k+2},\ldots , e_n+\varepsilon_n],\\ \hskip1cm\text{ if }\ r=\sum_{j=k+1}^n\varepsilon_jp^{e_j}
 \text{ with } \varepsilon_j\in \{0,1\},\\
0,\hskip.6cm\text{ otherwise.}
\end{cases}$$
\end{prop}

This proposition can easily be proved by using the Laplace development, the Cartan formula and the relations \ref{(i)} and \ref{(ii)}.

From the formula in \fullref{thm3.1}, one gets Wilkerson's formula as follows.

\begin{thm}{\rm (Wilkerson \cite{14})}\qua\label{thm3.4} For any integers $0\le s < n \le i,$ we have
$$\St^{\Delta_{i+1}}Q_{n,s}= P^{p^i}\St^{\Delta_i}Q_{n,s}.$$
\end{thm}

\begin{proof} Applying \fullref{thm3.1}, the Cartan formula and \fullref{prop3.3}, we get
\begin{align*}
P^{p^i}\St^{\Delta_i}Q_{n,s} &= (-1)^nP^{p^i}([0,1,\ldots,\hat s,\ldots, n-1,i]L_n^{p-2})\\
&= (-1)^n\sum_rP^r([0,1,\ldots,\hat s,\ldots, n-1,i])P^{p^i-r}(L_n^{p-2}),
\end{align*}
where the sum runs over all $$r=\varepsilon_0p^0 +\varepsilon_1p^1 +\ldots+\varepsilon_{s-1}p^{s-1}+\varepsilon_{s+1}p^{s+1}+\ldots + \varepsilon_{n-1}p^{n-1}+\varepsilon_ip^i$$ with $\varepsilon_j\in\{0,1\}$ for any $j$ and $r\le p^i$.

If $\varepsilon_i=0$ then $r < p^0+p^1+\ldots + p^{n-1}$ and 
\begin{align*}
2(p^i-r)& > 2(p^i-(p^0+p^1+\ldots + p^{n-1}))\\
&= 2(p^i-p^n+1+(p-2)(p^0+p^1+\ldots + p^{n-1}))\\
&> 2(p-2)(p^0+p^1+\ldots + p^{n-1}) = \dim L_n^{p-2}.
\end{align*}
This implies $P^{p^i-r}(L_n^{p-2})=0$. 

Since $r\le p^i$, if  $\varepsilon_i=1$ then $\varepsilon_j=0$ for $j\ne i$ and $r=p^i$. Hence, using the above equalities and \fullref{prop3.3}, we obtain
$$\eqalignbot{
P^{p^i}\St^{\Delta_i}Q_{n,s}&= (-1)^nP^{p^i}([0,1,\ldots,\hat s,\ldots, n-1,i])L_n^{p-2}\cr
&=(-1)^n [0,1,\ldots,\hat s,\ldots, n-1,i+1]L_n^{p-2}\cr
&= \St^{\Delta_{i+1}}Q_{n,s}.}
\proved$$
\end{proof}

Next, we compute the action of $\St^{\Delta_i}$ on M\`ui invariants.

\begin{thm}\label{thm3.5} For any positive integers $i,\ \! n,$ we have
$$\St^{\Delta_i}V_n = (-1)^{n-1}[0,1,\ldots,n-2,i]L_{n-1}^{p-2}.$$
\end{thm}

\begin{proof}  Since $L_n=L_{n-1}V_n$, applying the Cartan formula, we get
\begin{equation}\label{ct11}\St^{\Delta_i}L_{n}= L_{n-1}\St^{\Delta_i}V_{n}+V_{n}\St^{\Delta_i}L_{n-1}.
\end{equation}
Using \fullref{lem2.3} and \fullref{thm2.4}, we have
\begin{align*} 
\St^{\Delta_i}L_{n-1}&=[i,1,2,\ldots,n-2],\\
\St^{\Delta_i}L_n&=[i,1,2,\ldots,n-2,n-1]\\
&=\sum_{s=0}^{n-2}(-1)^{n+s}[i,1,2,\ldots,n-2,s]Q_{n-1,s}+[i,1,2,\ldots,n-2]V_n\\
&=(-1)^n[i,1,2,\ldots,n-2,0]Q_{n-1,0} +[i,1,2,\ldots,n-2]V_n.
\end{align*} 
Combining \eqref{ct11}, the above equalities and the relation $Q_{n-1,0}=L_{n-1}^{p-1}$, we get
$$\eqalignbot{\St^{\Delta_i}V_{n}&=  (\St^{\Delta_i}L_{n}-V_{n}\St^{\Delta_i}L_{n-1})/L_{n-1}\cr
&= (-1)^n[i,1,2,\ldots,n-2,0]Q_{n-1,0}/L_{n-1}\cr
&= (-1)^{n-1}[0,1,2,\ldots,n-2,i]L_{n-1}^{p-2}.}
\proved$$
\end{proof}

\begin{corl} For any  integers $0< i \le n,$ we have
$$\St^{\Delta_i}V_{n}=\begin{cases} 0,&i<n-1,\\
(-1)^{n-1}Q_{n-1,0}V_n,& i=n-1,\\
(-1)^{n-1}Q_{n-1,0}(Q_{n-1,n-2}^pV_n+V_n^p),&i=n.
\end{cases}$$
\end{corl}

\begin{proof}
If $i<n-1$ then $[0,1,\ldots,n-2,i]=0$.  Hence, $\St^{\Delta_i}V_{n}=0$. 

For $i=n-1$, we have $[0,1,2,\ldots,n-2,n-1]=L_n= L_{n-1}V_n$. Hence,  from \fullref{thm3.5}, we get
$$\St^{\Delta_{n-1}}V_{n}=(-1)^{n-1}L_{n-1}^{p-1}V_n= (-1)^{n-1}Q_{n-1,0}V_n.$$
Let $i=n$. A direct computation  shows 
\begin{align*}  [0,1,\ldots , n-2, n]&= L_{n,n-1}
= L_nQ_{n,n-1}\\
&=L_{n-1}V_n(Q_{n-1,n-2}^p + V_n^{p-1}). 
 \end{align*} 
From the above equalities, \fullref{thm3.5} and the relation $L_{n-1}^{p-1}=Q_{n-1,0}$, we obtain
$$\St^{\Delta_n}V_{n}=(-1)^{n-1}Q_{n-1,0}(Q_{n-1,n-2}^pV_n+V_n^p).$$
The corollary follows.
\end{proof}

\begin{thm} Set $s_0 = 0$. Then $\St^{\Delta_i}M_{n,s_1,\ldots,s_k}^{(d)}$ equals
$$ \begin{cases} (-1)^{s_t-t}M_{n,s_0,\ldots,\hat s_t,\ldots,s_k}^{(d)},&s_1 > 0, i = s_t, 1\le t\le k,\\
(-1)^{n-1}(d-1)M_{n,s_1,\ldots,s_k}[1,2,\ldots,n-1,i]L_n^{d-2},
&i\ge n, s_1=0,\\
(-1)^{n-1}\big((-1)^k[k;1,\ldots,\hat s_1,\ldots,\hat s_k,\ldots,n-1,i]L_n^{d-1}\\
\quad +(d-1)M_{n,s_1,\ldots,s_k}[1,2,\ldots,n-1,i]L_n^{d-2}\big),&
i\ge n, s_1>0,\\
           0,&\text{otherwise.}
            \end{cases}$$
\end{thm}

\begin{proof} Applying \fullref{lem2.2}, we have
$$\St^{\Delta_i}M_{n,s_1,\ldots,s_k} =\begin{cases} [k;i,1,\ldots,\hat s_1,\ldots,\hat s_k,\ldots,n-1],& s_1>0,\\ 0,&s_1=0.\end{cases}$$
If $i= s_t$ then
$[k;i,1,\ldots,\hat s_1,\ldots,\hat s_k,\ldots,n-1]= (-1)^{s_t-t}M_{n,s_0,\ldots,\hat s_t,  \ldots,s_k}.$

If $i\ge n$ then
$$[k;i,1,\ldots,\hat s_1,\ldots,\hat s_k,\ldots,n-1]=(-1)^{n-k-1}[k;1,\ldots,\hat s_1,\ldots, \hat s_k,\ldots,n-1,i].$$
Thus the theorem is proved for $d=1$.

 For $d>1$, using \fullref{lem2.2} and the Cartan formula, we have 
\begin{align*} \St^{\Delta_i}L_n^{d-1}&= (d-1)L_n^{d-2}\St^{\Delta_i}L_n,\\  \St^{\Delta_i}L_n&=(-1)^{n-1}[1,2,\ldots,n-1,i],\\
\St^{\Delta_i}M_{n,s_1,\ldots,s_k}^{(d)}&=\St^{\Delta_i}(M_{n,s_1,\ldots,s_k})L_n^{d-1}
+(d-1)M_{n,s_1,\ldots,s_k}L_n^{d-2}\St^{\Delta_i}L_n.
\end{align*}
Combining the above equalities we obtain the theorem. 
\end{proof}

\begin{thm}  For $1\le d\le p-1$, we have
$$ \St_uM_{n,s_1,\ldots,s_k}^{(d)} =\begin{cases} (-1)^{k+s_t-t}M_{
n,s_1,\ldots,\hat s_t,\ldots,s_k}^{(d)},&u = s_t,\\
(-1)^{n-1}[k-1;0,\ldots,\hat s_1,\ldots,\hat s_k,\ldots,n-1,u]L_n^{d-1},&
u\ge n,\\
           0,&\text{otherwise.} \end{cases}$$
\end{thm}

\begin{proof} Since $M_{n,s_1,\ldots,s_k}=[k;0,\ldots,\hat s_1,
\ldots,\hat s_k,\ldots,n-1]$, applying \fullref{lem2.2}, we get
$$\St_uM_{n,s_1,\ldots,s_k}=(-1)^{k-1}[k-1;u,0,\ldots,\hat s_1,\ldots,\hat s_k,\ldots,n-1].$$
If $0\le u \le n-1$ then
$$[k-1;u,0,\ldots,\hat s_1,\ldots,\hat s_k,\ldots,n-1]=\begin{cases}
(-1)^{s_t-t+1}M_{n,s_1,\ldots,\hat s_t,\ldots,s_k},&u = s_t,\\
0, &\text{otherwise}.\end{cases}$$
If  $u>n-1$ then we have
\begin{multline*} [k-1;u,0,\ldots,\hat s_1,\ldots,\hat s_k,\ldots,n-1]\\
= (-1)^{n-k}[k-1;0,\ldots,\hat s_1,\ldots,\hat s_k,\ldots,n-1,u].
\end{multline*}
The theorem is proved for $d=1$. 

Since $\St_uL_n=0$, using the Cartan formula, we get
$$\St_u(M_{n,s_1,\ldots,s_k}^{(d)})=\St_u(M_{n,s_1,\ldots,s_k})L_n^{d-1}.$$
The theorem now follows from the above equalities.  
\end{proof}

By the analogous argument as given in the proof of \fullref{thm3.4}, we can show that the Wilkerson formula also holds for M\`ui invariants.

\begin{thm} For any integers $i, u\ge n$, we have
\begin{align*}
\St^{\Delta_{i}}V_n&= P^{p^{i-1}}\St^{\Delta_{i-1}}V_n,\\
\St^{\Delta_{i+1}}M_{n,s_1,\ldots,s_k}^{(d)}&= P^{p^i}\St^{\Delta_i}M_{n,s_1,\ldots,s_k}^{(d)},\\
\St_{u+1}M_{n,s_1,\ldots,s_k}^{(d)} &= P^{p^u}\St_uM_{n,s_1,\ldots,s_k}^{(d)}.
\end{align*}
\end{thm}

\begin{rem} 
Using \fullref{thm2.4} and the above results, we can  compute the action of the primitive Steenrod--Milnor operations on the modular invariants  in terms of Dickson and M\`ui invariants for $i, u\ge n$. For example, by a direct calculation, we easily obtain
\begin{align*} 
\St^{\Delta_{n+1}}Q_{n,s} &= (-1)^nQ_{n,0}(Q_{n,n-1}^{p}Q_{n,s}-Q_{n,s-1}^p),\\
\St^{\Delta_{n+2}}Q_{n,s} &= (-1)^nQ_{n,0}(Q_{n,n-1}^{p^2+p}Q_{n,s}- Q_{n,n-2}^{p^2}Q_{n,s} + Q_{n,s-2}^{p^2} - Q_{n,s-1}^pQ_{n,n-1}^{p^2}).
\end{align*} 
Here, by convention, $Q_{n,t}=0$ for $t<0$.
\begin{align*} 
\St^{\Delta_{n+1}}V_n  = (-&1)^{n-1}Q_{n-1,0}\big((Q_{n-1,n-2}^{p^2+p} - Q_{n-1,n-3}^{p^2})V_n + Q_{n-1,n-2}^{p^2}V_n^p + V_n^{p^2}\big),\\
\St_nM_{n,s_1,\ldots,s_k}^{(d)} &=\sum_{t=1}^k(-1)^{n-1+k-t} M_{n,s_1,\ldots,\hat s_t,\ldots,s_k}^{(d)}Q_{n,s_t},\\
\St^{\Delta_n}M_{n,s_1,\ldots,s_k}^{(d)} &=(-1)^{n-1}\Big(\sum_{t=1}^k(-1)^{t} M_{n,s_0,\ldots,\hat s_t,\ldots,s_k}^{(d)}Q_{n,s_t} + dM_{n,s_1,\ldots,s_k}^{(d)}Q_{n,0}\Big),
\end{align*}
where $s_0=0$ and $s_1>0$. If $s_1=0$ then
$$\St^{\Delta_n}M_{n,s_1,\ldots,s_k}^{(d)} =(-1)^{n-1}(d-1)M_{n,s_1,\ldots,s_k}^{(d)}Q_{n,0}.$$
Furthermore, the computation of the action of the primitive Steenrod--Milnor operations on the modular invariants in terms of Dickson and M\`ui invariants by the use of our results in this section is more convenient than that by using Wilkerson's formula. For example, to compute $\St^{\Delta_{n+2}}Q_{n,s}$  by using Wilkerson's formula, we need to compute  $\smash{P^{\!p^{n+1}}(Q_{n,0}(Q_{n,n-1}^{p}Q_{n,s}-Q_{n,s-1}^p))}$ in terms of Dickson invariants.
But computing $\smash{P^{\!p^{n+1}}(Q_{n,0}(Q_{n,n-1}^{p}Q_{n,s}-Q_{n,s-1}^p))}$  is more difficult than that of $[0,1,\ldots,\hat s,\ldots,n-1,n+2]$.
\end{rem}

\section{On the description of the determinant invariants in terms of Dickson and M\`ui invariants}
\label{sec4}
In this section, we study the problem of description of the determinant invariants in terms of Dickson and M\`ui invariants. The explicit formulae for the determinant invariants in terms of Dickson and M\`ui invariants are useful tools for computing the action of the cohomology operations on the modular invariants.

In general, it is difficult to give  explicit formulae for this problem. In particular, for $n=2,3$, we can explicitly compute $[u,v], [u,v,w]$ in terms of M\`ui invariants and $[u,v], [u,v,v+1]$ in terms of Dickson invariants, where $u,v,w$ are nonnegative integers.

Note that the problem of description of $[u,v,w]$ in terms of Dickson invariants is complicated. It is still open.

\begin{prop} For $0\le u < v < w$, we have
\begin{align}\label{ct7}[u,v] &= \sum_{s=u}^{v-1}V_1^{p^v-p^{s+1}+p^u} V_2^{p^s}, \\
\label{ct8}{[u,v,w]} &= \sum_{s=u}^{v-1}[u,s+1][v,w]L_2^{-p^{s+1}}V_3^{p^s} + \sum_{s=v}^{w-1}[u,v][s+1,w]L_2^{-p^{s+1}}V_3^{p^s}.
\end{align}
\end{prop}

\begin{proof} The relation \eqref{ct7} is proved by induction on $v$. We prove \eqref{ct8} by induction on $v,w$. Applying \fullref{thm2.4}, we can easily prove the following by induction on $v$
\begin{equation}\label{ct9}[u,v,v+1] = \sum_{s=u}^{v-1}[u,s+1]L_2^{p^v-p^{s+1}}V_3^{p^s}.\end{equation}
Since $L_2^{p^v} = [v,v+1]$, the relation \eqref{ct8} holds for $w=v+1$.

Let $w=v+2$. By a direct computation using \fullref{thm2.4} and \eqref{ct9}, we have
\begin{align*}[u,v,v+2] &= [u,v,v+1]Q_{2,1}^{p^v} + [u,v]V_3^{p^v}\\
&= \sum_{s=u}^{v-1}[u,s+1]L_2^{p^v-p^{s+1}}V_3^{p^s}Q_{2,1}^{p^v} + [u,v]V_3^{p^v}. 
\end{align*}
We observe that $(L_2Q_{2,1})^{p^v} = [v,v+2]$, $L_2^{p^{v+1}}=[v+1,v+2]$. Hence, the relation \eqref{ct8} holds for $w=v+2$.
\allowdisplaybreaks
Suppose that \eqref{ct8} holds for $w$ and $w+1$. It is easy to see that $$[w+1,w]Q_{2,0}^{p^w} = - L_2^{p^{w+1}}.$$ Hence, using \fullref{thm2.4} and the inductive hypothesis, we get
\begin{align*}
[u,v,w+2] &= [u,v,w+1]Q_{2,1}^{p^w} - [u,v,w]Q_{2,0}^{p^w}+ [u,v]V_3^{p^w}\\
&= \Big(\sum_{s=u}^{v-1}[u,s+1][v,w+1]L_2^{-p^{s+1}}V_3^{p^s}\\
&\qquad+ \sum_{s=v}^{w}[u,v][s+1,w+1]L_2^{-p^{s+1}}V_3^{p^s}\Big)Q_{2,1}^{p^w}\\
&\qquad-\Big(\sum_{s=u}^{v-1}[u,s+1][v,w]L_2^{-p^{s+1}}V_3^{p^s}\\
&\qquad+\sum_{s=v}^{w-1}[u,v][s+1,w]L_2^{-p^{s+1}}V_3^{p^s}\Big)Q_{2,0}^{p^w} + [u,v]V_3^{p^w}
\\
&= \sum_{s=u}^{v-1}[u,s+1]\big([v,w+1]Q_{2,1}^{p^w} - [v,w]Q_{2,0}^{p^w}\big)L_2^{-p^{s+1}}V_3^{p^s}\\
&\qquad+\sum_{s=v}^{w}[u,v]\big([s+1,w+1]Q_{2,1}^{p^w} - [s+1,w]Q_{2,0}^{p^w}\big)L_2^{-p^{s+1}}V_3^{p^s}.
\end{align*}
This equality and  \fullref{thm2.4} imply the relation \eqref{ct8} for $w+2$, completing the proof.
\end{proof}

Now, we compute $[u,v]$ in terms of $L_2$ and $Q_{2,1}$.

Let $\alpha_i(a)$ denote the $i$--th coefficient in $p$--adic expansion of a nonnegative integer $a$. That means
$$a= \alpha_0(a)p^0+\alpha_1(a)p^1+\alpha_2(a)p^2+ \ldots ,$$ 
for $0 \le \alpha_i(a) <p, i\ge 0.$ We set $\alpha_i(a)=0$ for $i<0$. 

Denote by $I(u,v)$ the set of all integers $a$ satisfying
$$\begin{array}{ll}\alpha_i(a)+\alpha_{i+1}(a) \le 1,&\text{ for any }  i,\\
\alpha_i(a)=0, &\text{ for either } i<u \text{ or } i\ge v-2.
\end{array}$$
The following was proved in Sum \cite{9} for $p=2$.

\begin{prop} Under the above notation, we have
$$ [u,v] = \sum_{a\in I(u,v)}(-1)^aL_2^{p^u+p(p-1)a}Q_{2,1}^{\frac{p^{v-1}-p^u}{p-1}-(p+1)a}.$$
\end{prop}

\begin{proof} The proof is by induction on $v$. Obviously, $I(u,u+1)=I(u,u+2)=\{0\}$ and $[u,u+1] = L_2^{p^u}$, $[u,u+2] = L_2^{p^u}Q_{2,1}^{p^u}$. Hence, the proposition follows with $v=u+1$ and $v=u+2$. From the definition of the set $I(u,v)$, we obtain
\begin{equation} \label{ct10}I(u,v+2) = I(u,v+1)\cup (p^{v-1}+I(u,v)),\end{equation}
where $p^{v-1}+I(u,v) = \{p^{v-1}+a\ ;\ a \in I(u,v)\}$.

Combining \fullref{thm2.4}, the inductive hypothesis  and the relation $Q_{2,0} = L_2^{p-1}$, we get
\begin{align*}
[u,v+2] &= [u,v+1]Q_{2,1}^{p^v} - [u,v]Q_{2,0}^{p^v}\\
&= \Big(\sum_{a\in I(u,v+1)}(-1)^aL_2^{p^u+p(p-1)a}Q_{2,1}^{\frac{p^{v}-p^u}{p-1}-(p+1)a}\Big)Q_{2,1}^{p^v}\\
&\quad-\Big(\sum_{a\in I(u,v)}(-1)^aL_2^{p^u+p(p-1)a}Q_{2,1}^{\frac{p^{v-1}-p^u}{p-1}-(p+1)a}\Big)Q_{2,0}^{p^v}\\
&= \sum_{a\in I(u,v+1)}(-1)^aL_2^{p^u+p(p-1)a}Q_{2,1}^{\frac{p^{v+1}-p^u}{p-1}-(p+1)a}\\
&\quad+\sum_{a\in I(u,v)}(-1)^{p^{v-1}+a}L_2^{p^u+p(p-1)(p^{v-1}+a)}Q_{2,1}^{\frac{p^{v+1}-p^u}{p-1}-(p+1)(p^{v-1}+a)}.
\end{align*}
From this equality and \eqref{ct10}, we see that the proposition is true for $v+2$, so the proof is completed.
\end{proof}

Now, we compute $[u,v,v+1]$ in terms of $L_3, Q_{3,1}, Q_{3,2}$.

Denote by $J(u,v)$ the set of all integers $a$ satisfying
$$\begin{array}{ll}\alpha_i(a) \le 1 \quad \text{and} \quad \alpha_i(a)+\alpha_{i+1}(a)+ \alpha_{i+2}(a) \le 2, &\text{ for any } i,\\
\alpha_i(a)=0, &\text{ for either } i<u \text{ or } i \ge v-2.
\end{array}$$
It is easy to see that for any $a\in J(u,v)$, there exists uniquely an expansion
$$a=a_0+p^{i_1}+p^{i_1+1}+a_1+\ldots +p^{i_k}+p^{i_k+1}+a_k,$$
with $i_0=u-3<i_1<\ldots <i_k<i_{k+1}=v-1, i_{j+1}-i_j\ge 3$ and $a_j \in I(i_j+3,i_{j+1})$ for $0\le j \le k$. 

We define the functions $b_{u,v}, c_{u,v}\co J(u,v) \to \mathbb Z$ by setting
\begin{align*} b_{u,v}(a) &= \frac{p^{v-1}-p^u}{p-1}-(p+1)a+p(p^{i_1}+ \ldots+p^{i_k}),\\
c_{u,v}(a)&=a_0+a_1+\ldots+a_k.
\end{align*}

\begin{prop}\label{prop4.3} Under the above notation, we have
$$ [u,v,v+1] = \sum_{a\in J(u,v)}(-1)^aL_3^{p^u+p(p-1)a}Q_{3,1}^{b_{u,v}(a)}Q_{3,2}^{c_{u,v}(a)}.$$
\end{prop}

The proof of the proposition is based on some lemmas.

\begin{lem}\label{lem4.4} For $0\le u<v$,
$$J(u,v+3) = J(u,v+2)\cup (p^v+J(u,v+1))\cup (p^v+p^{v-1}+J(u,v)).$$
Here, for $x \in \mathbb Z$ and $A\subset \mathbb Z$, we write $x+A=\{x+a \ ;\ a\in A\}$.
$$\begin{array}{rll}
b_{u,v+3}(a) =&p^{v+1}+b_{u,v+2}(a),& \cr
c_{u,v+3}(a) =&c_{u,v+2}(a),&\text{ for } a\in J(u,v+2),\cr
b_{u,v+3}(p^v+ a) =&b_{u,v+1}(a),& \cr
c_{u,v+3}(p^v+ a) =&p^v + c_{u,v+1}(a),& \text{ for  } a\in J(u,v+1),\cr
b_{u,v+3}(p^v+p^{v-1}+a) =&b_{u,v}(a),& \cr
c_{u,v+3}(p^v+p^{v-1}+a) =&c_{u,v}(a),&\text{ for  } a\in J(u,v).
\end{array}$$
\end{lem}
This lemma can easily be proved by computing directly from the definitions of $J(u,v)$, $b_{u,v}$ and $c_{u,v}$.

\begin{lem}\label{lem4.5} For any $0\le u<v$, we have
\begin{multline*}
[u,v+3,v+4] = [u,v+2,v+3]Q_{3,1}^{p^{v+1}} \\
-[u,v+1,v+2]Q_{3,0}^{p^{v+1}}Q_{3,2}^{p^v} + [u,v,v+1]Q_{3,0}^{p^{v+1}+p^v}.
\end{multline*}
\end{lem}

\begin{proof}  A direct calculation using \fullref{thm2.4} gives
\begin{align*}
[u,v+3,v+4] &= [u,v+2,v+3]Q_{2,0}^{p^{v+2}} + [u,v+3]V_3^{p^{v+2}}\\
&= [u,v+2,v+3](Q_{3,1}^{p^{v+1}} - Q_{2,1}^{p^{v+1}}V_3^{(p-1)p^{v+1}})\\
&\quad+([u,v+2]Q_{2,1}^{p^{v+1}} - [u,v+1]Q_{2,0}^{p^{v+1}})V_3^{p^{v+2}}\\
&\hskip3cm (\text{since } Q_{3,1}=Q_{2,0}^p+Q_{2,1}V_3^{p-1})\\
&=[u,v+2,v+3]Q_{3,1}^{p^{v+1}}\\
&\quad- ([u,v+1,v+2]Q_{2,0}^{p^{v+1}}+ [u,v+2]V_3^{p^{v+ 1}})Q_{2,1}^{p^{v+1}}V_3^{(p-1)p^{v+1}}\\
&\quad+ [u,v+2]Q_{2,1}^{p^{v+1}}V_3^{p^{v+2}} - [u,v+1]Q_{2,0}^{p^{v+1}}V_3^{p^{v+2}}\\
&=[u,v+2,v+3]Q_{3,1}^{p^{v+1}}\\
&\quad-[u,v+1,v+2]Q_{2,0}^{p^{v+1}}V_3^{(p-1)p^{v+1}}(Q_{2,1}^{p^{v+1}}+V_3^{(p-1)p^v})\\
&\quad+([u,v+1,v+2] - [u,v+1]V_3^{p^v})Q_{2,0}^{p^{v+1}}V_3^{(p-1)(p^{v+1}+p^v)}.
\end{align*}
Using \fullref{thm2.4} and the relations $Q_{3,2}=Q_{2,1}^p+V_3^{p-1}$, $Q_{3,0}=Q_{2,0}V_3^{p-1}$, we obtain the lemma.
\end{proof}

\begin{proof}[Proof of \fullref{prop4.3}] The proof is by induction on $v$. For $v=u+1, u+2, u+3$ the proposition is obvious. Suppose that it is true for $v, v+1, v+2.$ 
 Using \fullref{lem4.5},  the inductive hypothesis and the relation $Q_{3,0} = L_3^{p-1}$, we get
\begin{align*}
[u,v+3,v+4] &= \sum_{a\in J(u,v+2)}\!\!(-1)^aL_3^{p^u+p(p-1)a}Q_{3,1}^{p^{v+1}+b_{u,v+2}(a)}Q_{3,2}^{c_{u,v+2}(a)}\\
&\   +\!\!\!\sum_{a\in J(u,v+1)}\!(-1)^{p^v+a}L_3^{p^u+p(p-1)(p^v+a)}Q_{3,1}^{b_{u,v+1}(a)}Q_{3,2}^{p^v+c_{u,v+1}(a)}\\
&\   +\!\!\!\sum_{a\in J(u,v)}\!(-1)^{p^v+p^{v-1}+a}L_3^{p^u+p(p-1)(p^v+p^{v-1}+a)}Q_{3,1}^{b_{u,v}(a)}Q_{3,2}^{c_{u,v}(a)}.
\end{align*}
Combining this equality and \fullref{lem4.4}, we see that the proposition holds for $v+3$. Hence, the proposition is proved. 
\end{proof}

\bibliographystyle{gtart}
\bibliography{link}

\begin{thebibliography}{}
\providecommand\bibmarginpar{\leavevmode\marginpar}
\def\urlstyle#1{{\tt #1}}

\bibitem{1}
\textbf{L\,E Dickson},
  \href{http://links.jstor.org/sici?sici=0002-9947(191101)12:1%3C75:AFSOIO%3E2%
.0.CO%3B2-%23} {\emph{A fundamental system of invariants of the general modular
  linear group with a solution of the form problem}}, Trans. Amer. Math. Soc.
  12 (1911) 75--98 \xox{MR}{1500882} \xox{JFM}{42.0136.01}

\bibitem{5}
\textbf{N\,H\,V H\uhorn{}ng}, \textbf{P\,A Minh}, \emph{The action of the mod
  {$p$} {S}teenrod operations on the modular invariants of linear groups},
  Vietnam J. Math. 23 (1995) 39--56 \xox{MR}{1367491}

\bibitem{4}
\textbf{N\,E Kechagias},
  \href{http://links.jstor.org/sici?sici=0002-9939(199307)118:3%3C943:TSAAOG%3%
E2.0.CO%3B2-K} {\emph{The {S}teenrod algebra action on generators of rings of
  invariants of subgroups of $\mathrm{GL}_n(\mathbb{Z}/p\mathbb{Z})$}}, Proc.
  Amer. Math. Soc. 118 (1993) 943--952 \xox{MR}{1152986}

\bibitem{6}
\textbf{I Madsen}, \textbf{R\,J Milgram}, \emph{The classifying spaces for
  surgery and cobordism of manifolds}, Annals of Mathematics Studies 92,
  Princeton University Press, Princeton, N.J. (1979) \xox{MR}{548575}

\bibitem{7}
\textbf{J Milnor},
  \href{http://links.jstor.org/sici?sici=0003-486X(195801)2:67:1%3C150:TSAAID%%
3E2.0.CO%3B2-2} {\emph{The {S}teenrod algebra and its dual}}, Ann. of Math.
  $(2)$ 67 (1958) 150--171 \xox{MR}{0099653}

\bibitem{2}
\textbf{H M\`ui}, \emph{Modular invariant theory and cohomology algebras of
  symmetric groups}, J. Fac. Sci. Univ. Tokyo Sect. IA Math. 22 (1975) 319--369
  \xox{MR}{0422451}

\bibitem{3}
\textbf{H M\`ui}, \href{http://dx.doi.org/10.1007/BF01163361} {\emph{Cohomology
  operations derived from modular invariants}}, Math. Z. 193 (1986) 151--163
  \xox{MR}{852916}

\bibitem{8}
\textbf{M\,D Neusel}, \emph{Inverse invariant theory and {S}teenrod
  operations}, Mem. Amer. Math. Soc. 146 (2000) x+158 \xox{MR}{1693799}

\bibitem{12}
\textbf{L Smith}, \textbf{R\,M Switzer},
  \href{http://links.jstor.org/sici?sici=0002-9939(198310)89:2%3C303:RANODA%3E%
2.0.CO%3B2-K} {\emph{Realizability and nonrealizability of {D}ickson algebras
  as cohomology rings}}, Proc. Amer. Math. Soc. 89 (1983) 303--313
  \xox{MR}{712642}

\bibitem{13}
\textbf{N\,E Steenrod}, \textbf{D\,B\,A Epstein}, \emph{Cohomology operations},
  Annals of Mathematics Studies 50, Princeton University Press, Princeton, N.J.
  (1962) \xox{MR}{0145525}

\bibitem{9}
\textbf{N Sum}, \emph{On the action of the {S}teenrod--{M}ilnor operations on
  the modular invariants of linear groups}, Japan. J. Math. $($N.S.$)$ 18
  (1992) 115--137 \xox{MR}{1173832}

\bibitem{10}
\textbf{N Sum}, \emph{On the action of the {S}teenrod algebra on the modular
  invariants of special linear group}, Acta Math. Vietnam. 18 (1993) 203--213
  \xox{MR}{1292080}

\bibitem{11}
\textbf{N Sum},
  \href{http://projecteuclid.org/getRecord?id=euclid.kmj/1138040053}
  {\emph{Steenrod operations on the modular invariants}}, Kodai Math. J. 17
  (1994) 585--595 \xox{MR}{1296929}Workshop on Geometry and Topology (Hanoi,
  1993)

\bibitem{14}
\textbf{C Wilkerson}, \emph{A primer on the {D}ickson invariants}, from:
  ``Proceedings of the Northwestern Homotopy Theory Conference (Evanston, Ill.,
  1982)'', Contemp. Math. 19, Amer. Math. Soc., Providence, RI (1983)  421--434
  \xox{MR}{711066}

\end{thebibliography}

\end{document}